\documentclass[11pt]{article}

\usepackage{amsmath,amssymb,amsthm,mathtools}
\usepackage{geometry}
\usepackage{enumitem}
\usepackage[hidelinks]{hyperref}
\usepackage{microtype}
\usepackage{mathrsfs}
\numberwithin{equation}{section}

\newtheorem{theorem}{Theorem}[section]
\newtheorem{lemma}{Lemma}[section]

\newtheorem{remark}{Remark}[section]

\newcommand{\R}{\mathbb{R}}
\newcommand{\Rp}{\mathbb{R}_{+}}

\newcommand{\eps}{\varepsilon}
\newcommand{\cR}{\mathcal{R}}

\newcommand{\dist}{\operatorname{dist}}
\newcommand{\norm}[1]{\left\|#1\right\|}

\title{Global dynamics of impulsive reaction--diffusion systems\\
with asymptotically unbounded domain}
\author{Xiao-Qiang Zhao\\ 
	Department of Mathematics and Statistics\\
	Memorial University of Newfoundland\\
	St. John's, NL A1C  5S7, Canada\\
	E-mail:\, zhao@mun.ca}

\date{}

\begin{document}
\maketitle

\begin{abstract}
This paper is devoted to the study of the global dynamics of  continuous--discrete hybrid population models with an asymptotically unbounded domain.  We first prove a global attractivity result for a class of impulsive reaction--diffusion equations whose diffusion coefficients vanish asymptotically. A cut-off function is combined with an auxiliary impulsive kinetic system to construct a subsolution that is compatible with the birth update at every cycle. We then use comparison arguments to establish threshold-type dynamics for the original model when the birth function is either monotone or nonmonotone. It turns out that the threshold value for extinction and persistence is explicitly determined by the birth pulse, the within-cycle growth, and the asymptotic expansion rate of the habitat.
\end{abstract}

\medskip
\noindent\textbf{Mathematical Subject Classification.} 35K57, 37N25, 92D25.

\medskip
\noindent\textbf{Keywords.} Impulsive systems, time-varying domains, asymptotically vanishing diffusion, threshold dynamics.

\section{Introduction}

A seasonal population may undergo qualitatively different mechanisms during different parts of one life cycle. Reproduction can be concentrated in a short interval, whereas movement, mortality, and competition act over the remainder of the cycle. This separation of time scales leads naturally to hybrid models in which a discrete update is coupled with a continuous spatial evolution. One of the basic formulations of this type was introduced by Lewis and Li \cite{LewisLi2012}, and it was extended to the cases of  higher-dimensional diffusion and advection \cite{FazlyLewisWang2017,FazlyLewisWang2020}. The same continuous--discrete viewpoint has  been used for populations with nonlocal movement \cite{WuZhao2019}, heterogeneous landscapes \cite{WuZhao2022}, synchronized emergence \cite{BaiLouZhao2022}, and climate-driven habitat shifts \cite{ZhangYiChen2024}. An important mathematical question in these models is how the discrete birth mechanism and the continuous dispersal stage combine to determine extinction, persistence, and spatial spread.

A second source of nonautonomy arises when the region available to the population is itself changing. The geometry of a habitat may vary periodically, approach a terminal size, or continue to expand over long time intervals, see, e.g., \cite{LamZhaoZhu2024} and references therein.
When such geometric change is prescribed by environmental conditions rather than by the population density, it is natural to regard the domain evolution as an external component of the model. Impulsive reaction--diffusion equations on evolving domains have already been considered in several settings. For example, periodic domain change with harvesting was studied in \cite{MengLinPedersen2021}, while a logistic model with harvesting pulses on a growing domain was invdestigated in
\cite{ZhangZhu2025}. More recently, general birth-pulse systems were treated for asymptotically bounded domains in \cite{ZhuZhao2025} and asymptotically periodic domains in \cite{FanZhao2026}, respectively.
However, it remains an unsolved problem to characterize the global dynamics of such impulsive systems with an asymptotically unbounded domain.

The asymptotically unbounded domain studied in this paper  is an expanding habitat whose size tends to infinity and whose relative expansion rate has a finite limit. After rescaling such a domain to a fixed one, the diffusion coefficient vanishes asymptotically, while the dilution term may converge to a nonzero constant. Consequently, the evolution system under consideration  is not a perturbation of a limiting parabolic equation since its spatial part degenerates and the relevant limiting object is an impulsive kinetic system.  Lam, Zhao and Zhu \cite{LamZhaoZhu2024} treated the analogous vanishing-diffusion mechanism for nonimpulsive reaction--diffusion systems and obtained interior convergence by specially constructed subsolutions. To the best of our knowledge, a corresponding  result for impulsive reaction--diffusion equations with asymptotically vanishing diffusion, covering time-periodic KPP-type reactions and birth-pulse maps, has not been established. The main analytical contribution of the present paper is to obtain such a result and then apply it to  implusive reaction-diffusion models on  an  asymptotically unbounded domain. In particular, the  subsolution must be compatible not only with the parabolic evolution but also with the birth update at every cycle.

To formulate   impulsive reaction-diffusion models on  time-varying domains, we write real time \(\hat t\ge0\)  uniquely as
\begin{equation*}
	\hat t=mT+t,
	\qquad
	m=\left\lceil\frac{\hat t}{T}\right\rceil-1,
	\qquad
	t\in(0,T],
\end{equation*}
where \(T>0\) is the length of one time cycle, and $\left\lceil\frac{\hat t}{T}\right\rceil$
denotes the smallest integer greater than or equal to $\frac{\hat t}{T}$.
Let \(H_m(x)\) be the population density at the beginning of the \(m\)-th cycle, and let \(v_m(t,x)\) be the density within that cycle. 
Following \cite{ZhuZhao2025,FanZhao2026}, we consider the following evolution system on \(\Omega(\hat t)\):
{\small
	\begin{equation}\label{model}
		\left\{
		\begin{aligned}
			\frac{\partial v_m}{\partial t}
			+\nabla v_m\cdot a(mT+t,x)
			+v_m\nabla\cdot a(mT+t,x)
			&=d\Delta v_m+f(mT+t,v_m),
			&&0<t\le T,\ x\in\Omega(mT+t),\\
			v_m(0,x)&=g(H_m(x)),
			&&x\in\Omega(mT),\\
			H_{m+1}(x)&=v_m(T,x),
			&&x\in\Omega((m+1)T),\\
			Bv_m(t,x)&=0,
			&&0<t\le T,\ x\in\partial\Omega(mT+t).
		\end{aligned}
		\right.
\end{equation}}Here \(d>0\) is the diffusion coefficient, \(f\) is the reaction term, \(g\) describes the birth pulse, and the operator \(B\) denotes either the homogeneous Dirichlet boundary condition or the homogeneous Neumann boundary condition. The vector field \(a\) represents the velocity generated by the movement of the habitat.

For simplicity, we consider a class of domains that evolve by linear isotropic deformation:
\begin{equation*}
	\Omega(\hat t)=\rho(\hat t)\Omega_0
	:=\{\rho(\hat t)y:y\in\Omega_0\},
\end{equation*}
where \(\Omega_0\subset\R^n\) is a bounded domain with \(C^{2+\theta}\)-boundary for some \(0<\theta\le1\), \(\rho\in C^1(\Rp,(0,\infty))\), and \(\rho(0)=1\). 
It then follows that 
\[
a(\hat t,x)=\frac{\rho'(\hat t)}{\rho(\hat t)}x,
\quad
\nabla\cdot a(\hat t,x)=n\frac{\rho'(\hat t)}{\rho(\hat t)}.
\]
Set
\begin{equation*}
	u_m(t,y):=v_m(t,\rho(mT+t)y),
	\quad
	N_m(y):=H_m(\rho(mT)y),
	\quad \forall  m\geq 0,\, y\in\Omega_0.
\end{equation*}
Thus, system \eqref{model} is reduced to
\begin{equation}\label{ReduceE}
	\left\{
	\begin{aligned}
		\frac{\partial u_m}{\partial t}
		&=\frac{d}{\rho^2(mT+t)}\Delta u_m
		+f(mT+t,u_m)
		-n\frac{\rho'(mT+t)}{\rho(mT+t)}u_m,
		&&0<t\le T,\ y\in\Omega_0,\\
		u_m(0,y)&=g(N_m(y)),
		&&y\in\Omega_0,\\
		N_{m+1}(y)&=u_m(T,y),
		&&y\in\Omega_0,\\
		Bu_m(t,y)&=0,
		&&0<t\le T,\ y\in\partial\Omega_0,
	\end{aligned}
	\right.
\end{equation}
where  $m\geq 0$, and $N_0(y)=\phi(y)\ge 0$.
As in \cite{FanZhao2026}, we impose the following assumptions on $f$ and $g$:
\begin{enumerate}
\item[(A1)] \(f\in C(\Rp\times\Rp,\R)\), \(f_u\in C(\Rp\times\Rp,\R)\), \(f(t+T,u)=f(t,u)\) for all \(t,u\ge0\), and \(f(t,0)\equiv 0\); for each  \(t\ge0\), \(f(t,u)/u\) is strictly decreasing for \(u>0\); there exists \(\check N>0\) such that \(f(t,N)<0\) for all \(N\ge\check N\) and \(t\in[0,T]\); and \(f\) is locally Lipschitz in \(u\), uniformly with respect to \(t\in[0,T]\).
\item[(A2)] 
\(g\in C(\Rp,\Rp)\), \(g(0)=0\), \(g'(0)>0\), and \(g(N)>0\) for all \(N>0\). The function \(g(N)/N\) is nonincreasing for \(N>0\), and \(g\) is locally Lipschitz on \(\Rp\).
\item[(A3)] 
Either \(g(\overline N)\le\overline N\) for some \(\overline N>0\), or 
$\limsup_{N\to\infty}\max_{t\in[0,T]}\frac{f(t,N)}{N^\nu}<0$ for some \(\nu>1\).
\item[(A4)] 
There exists \(\sigma>0\) such that \(g\) is nondecreasing on \([0,\sigma]\).
\end{enumerate}

Clearly, the time periodicity of $f$  implies that in system \eqref{ReduceE},  $f(mT+t,u_m)=f(t,u_m)$ for all $m\geq 0$.  For typical examples of $f$ and $g$, we refer to \cite{ZhuZhao2025,FanZhao2026}. 
Accordingly, we investigate the global dynamics of the evolution system \eqref{model} in an asymptotically unbounded domain. It suffices to study the impulsive system \eqref{ReduceE} under the assumption that the scaling factor $\rho(\hat t)$ satisfies
\begin{enumerate}
\item[(A5)] 
 $\rho  (0)=1$,  $\lim_{\hat t\to\infty}\rho(\hat t)=+\infty$,  and 
 $\lim_{\hat t\to\infty}\frac{\rho'(\hat t)}{\rho(\hat t)}=k\ge 0$.
\end{enumerate}

The main difficulty in our analysis is to prove persistence for system \eqref{ReduceE} in the interior of the fixed reference domain. For this purpose, we first consider a large class of impulsive reaction--diffusion equations with asymptotically vanishing diffusion. Let \(\phi_\eta\) be a smooth cut-off function satisfying \(\Delta\phi_\eta\ge-a_\eta\phi_\eta\), and let \((w_m(t),M_m)\) solve an auxiliary impulsive kinetic system with a small additional mortality term. We show that
$(\underline u_m(t,x),\underline N_m(x)):=\phi_\eta(x)(w_m(t),M_m)$
is an impulsive subsolution. Besides the parabolic inequality, the construction must satisfy the inequality 
at each birth pulse, which follows from  the subhomogeneity of \(g\). This cut-off--kinetic construction is tailored to the impulsive setting and differs from the connecting-orbit construction 
proposed  in \cite{LamZhaoZhu2024} for the nonimpulsive problem.

The rest of this paper is organized as follows. In Section~2, we present preliminary results on the uniform boundedness of solutions for system \eqref{ReduceE}, the threshold dynamics of the limiting impulsive kinetic system, and the existence of a suitable cut-off function. In Section~3, we establish the global dynamics for a class of impulsive reaction--diffusion equations with asymptotically vanishing diffusion by constructing a subsolution adapted to the impulsive structure. In Section~4, we apply these results to the original impulsive model on an asymptotically unbounded domain when the birth function is monotone or nonmonotone, respectively. A brief discussion concludes the paper.

\section{Preliminaries}

Let $C(\overline{\Omega_0}, \mathbb{R})$ be the Banach space of all continuous functions from $\overline{\Omega_0}$ to
 $\mathbb{R}$ equipped with the maximum norm $\|\cdot \|_\infty$, and $C_0(\overline{\Omega_0}, \mathbb{R})=
 \{\phi\in C(\overline{\Omega_0}, \mathbb{R}):\, \, \phi(x)=0, \, \forall x\in \partial \Omega_0 \}$.
  Throughout this paper, 
 we choose $X=C(\overline{\Omega_0}, \mathbb{R})$ and $X_+=C(\overline{\Omega_0}, \mathbb{R}_+)$ in the case of  homogeneous Neumann boundary condition;
 $X=C_0(\overline{\Omega_0}, \mathbb{R})$
 and $X_+=C_0(\overline{\Omega_0}, \mathbb{R}_+)$ in the case of homogeneous Dirichlet boundary condition.
	
By the same arguments as in \cite[Lemma 2.1]{FanZhao2026}, we have the following result on the 
uniform boundedness of solutions of system \eqref{ReduceE}.

\begin{lemma}\label{lem:bounded}
	Let \emph{(A1)--(A5)}  hold. Then solutions of system \eqref{ReduceE} are uniformly bounded in the sense that  for every \(r>0\), there exists \(C(r)>0\) such that for any $\phi\in X_+$ with  \(\norm{\phi}_\infty\le r\), there holds
	$\norm{u_m(t,\cdot,\phi)}_\infty\le C(r)$ for all $m\ge0$ and $t\in[0,T]$.
\end{lemma}

In view of (A5), the diffusion coefficient in system \eqref{ReduceE} tends to zero and the dilution coefficient converges to \(nk\). We are therefore led to the following limiting impulsive system:
\begin{equation}
\left\{
\begin{aligned}
\frac{dw_m}{dt}&=f(t,w_m)-nkw_m,
&&0<t\le T,\\
w_m(0)&=g(M_m),\\
M_{m+1}&=w_m(T),
\end{aligned}
\right.
\label{eq:limiting-kinetic}
\end{equation}
where $m\geq 0$ and \(M_0\ge0\). Let \(W_k(t)\) be the solution map at time \(t\) of the $T$-periodic equation:
\begin{equation}
\frac{dw}{dt}=f(t,w)-nkw,
\qquad t>0,
\label{eq:kinetic-flow-k}
\end{equation}
and let  \(W_k=W_k(T)\) be the Poincar\'e map of system \eqref{eq:kinetic-flow-k}.  Then \eqref{eq:limiting-kinetic} generates the discrete-time recursion
\begin{equation}
M_{m+1}=W_k\big(g(M_m)\big)=:\widehat W_k(M_m),
\qquad m\ge0.
\label{eq:kinetic-map-k}
\end{equation}

The linearization of system \eqref{eq:limiting-kinetic} at the zero solution is
\begin{equation}
\left\{
\begin{aligned}
(z_m)'(t)&=(f_u(t,0)-nk)z_m(t),
&&0<t\le T,\\
z_m(0)&=g'(0)Z_m,\\
Z_{m+1}&=z_m(T).
\end{aligned}
\right.
\label{eq:linearized-kinetic}
\end{equation}
 It is easy to see that the spectral radius of the cycle map of linear system \eqref{eq:linearized-kinetic}  is
\begin{equation*}
\cR_k
:=g'(0)\exp\left\{\int_0^T\big(f_u(s,0)-nk\big)\,ds\right\}.
\end{equation*}
Then we have  the following threshold type result  in terms of $\cR_k$.

\begin{lemma}\label{lem:kinetic}
Assume \emph{(A1)--(A3)} and suppose that \(g\) is nondecreasing on \(\Rp\). Then the following statements are valid for system \eqref{eq:kinetic-map-k}:
\begin{enumerate}[label=\rm(\roman*)]
\item If \(\cR_k\le1\), then \(0\) is globally asymptotically stable in \(\Rp\).
\item If \(\cR_k>1\), then \(\widehat W_k\) admits a unique positive fixed point \(M_k^*>0\), which is
globally asymptotically stable in \(\Rp\setminus \{0\}\).
\end{enumerate}
\end{lemma}

\begin{proof}
Clearly, the scalar equation \eqref{eq:kinetic-flow-k} admits the comparison principle. Then 
monotonicity of $g$ on \(\Rp\) implies  that  the map \(\widehat W_k=W_k\circ g\) is monotone (nondecreasing) on
$\mathbb{R}_+$. Since \(f(t,u)/u\) is strictly decreasing for $u>0$, it easily follows that  for each $t\geq 0$, 
the function $h_k(t,u):=f(t,u)-nku$  is strictly subhomogeneous on $\mathbb{R}_+$ in the sense that $h_k(t,\alpha u)>\alpha h_k(t,u)$ whenever 
\(u>0\) and \(\alpha\in(0,1)\). 
Let \(w(t):=W_k(t)u\). It is easy to verify that  \(\alpha w(t)\) is a strict subsolution of \eqref{eq:kinetic-flow-k} with initial value \(\alpha u\). 
Thus,  the  comparison principle implies $W_k(T)(\alpha u)>\alpha W_k(T)u$, that is, 
\(W_k\) is strongly subhomogeneous on  $\mathbb{R}_+$. Since \(g(u)/u\) is nonincreasing for \(u>0\),  we 
see that $g(u)$ is subhomogeneous on $\mathbb{R}_+$
in the sense that  $g(\alpha u)\ge \alpha g(u)$ for all  $u\ge0$ and $\alpha\in[0,1]$.
This, together with the monotonicity and strong subhomogeneity of \(W_k\), gives rise to 
\[
\widehat W_k(\alpha u)
=W_k\bigl(g(\alpha u)\bigr)
\ge W_k\bigl(\alpha g(u)\bigr)
>\alpha W_k\bigl(g(u)\bigr)
=\alpha\widehat W_k(u),\quad \forall  u>0,\,  \alpha\in(0,1),
\]
and hence,  \(\widehat W_k\) is strongly subhomogeneous on $\mathbb{R}_+$.  It follows from 
the comparison argument used in Lemma~\ref{lem:bounded} that all forward orbits for  \(\widehat W_k\) 
are bounded.  By  the differentiability of solutions with repect to the initial value, we further  obtain
\[
\widehat W_k'(0)
=g'(0)\exp\left\{\int_0^T\bigl(f_u(s,0)-nk\bigr)\,ds\right\}
=\cR_k.
\]
The threshold type conclusions now follow from \cite[Theorem 2.3.4]{Zhao2017}.
\end{proof}

In the case where \(\cR_k>1\), define
\begin{equation*}
w_k^*(t):=W_k(t)\big(g(M_k^*)\big),
\qquad 0\le t\le T.
\end{equation*}
Then \(w_k^*(T)=M_k^*\), while \(w_k^*(0)=g(M_k^*)\). In the impulsive sense, \(w_k^*(t)\) is the unique positive \(T\)-periodic solution associated with \eqref{eq:limiting-kinetic}.  As a consequence of Lemma~\ref{lem:kinetic}, we have the following result for 
system \eqref{eq:limiting-kinetic}.

\begin{lemma}\label{lem:kineticflow}
Under the assumptions of Lemma~\ref{lem:kinetic}, the following statements are valid:
\begin{enumerate}[label=\rm(\roman*)]
\item If \(\cR_k\le1\), then every nonnegative solution of \eqref{eq:limiting-kinetic} satisfies
$\lim_{m\to\infty}w_m(t)=0$ uniformly for \(t\in[0,T]\).
\item If \(\cR_k>1\), then every positive solution of \eqref{eq:limiting-kinetic} satisfies
$\lim_{m\to\infty}w_m(t)=w_k^*(t)$ uniformly for \(t\in[0,T]\).
\end{enumerate}
\end{lemma}

We need the following elementary observation in our construction of a subsolution for impulsive systems with asymptotically vanishing diffusion.
\begin{lemma}\label{lem:cutoff}
Let \(\Omega\subset\R^n\)  be a bounded domain and 
assume that $U\subset\Omega$ is open and $K\subset U$ is compact. Then there exists a smooth cut-off function 
$\phi\in C_c^\infty(\Omega)$ and a number \(a>0\) such that
$0\leq \phi\leq 1$ in $\Omega$, \, $\phi\equiv 1$ on  $K$,\,  $\phi=0$
in $\Omega\setminus U$,\,  and  $\Delta\phi\geq -a\phi$ in $\Omega$.
\end{lemma}

\begin{proof}
	Define
	\[
	\psi(x)
	=
	\begin{cases}
		c_n\exp\!\left(-\dfrac{1}{1-|x|^2}\right),
		& |x|<1,\\[1.2ex]
		0,
		& |x|\geq1,
	\end{cases}
	\]
	where $c_n>0$ is chosen such that $\int_{\mathbb{R}^n}\psi(x)\,\mathrm dx=1$. Then $\psi\in C_c^\infty(\mathbb{R}^n)$ and $\psi\geq0$. 
	For $r=|x|<1$, a straightforward computation gives
	\[
	\frac{\Delta\psi(x)}{\psi(x)}
	=
	\frac{4r^2}{(1-r^2)^4}
	-\frac{8r^2}{(1-r^2)^3}
	-\frac{2n}{(1-r^2)^2}.
	\]
	The right-hand side is continuous for $0\leq r<1$ and tends to $+\infty$ as $r\to1^-$, and hence,  it is bounded from below. 
	Thus, there exists a constant $C_0>0$ such that $\Delta\psi\geq-C_0\psi$ in $B_1(0)$. Since $\psi\equiv0$ outside $B_1(0)$ and $\psi\in C_c^\infty(\mathbb{R}^n)$, the same inequality holds in $\mathbb{R}^n$.

	Since $K\subset U$ is compact, we have $d:=\operatorname{dist}\bigl(K,\mathbb{R}^n\setminus U\bigr)>0$. Choose $\varepsilon>0$ sufficiently small so that $3\varepsilon<d$. For $r>0$, set $K_r:=\bigl\{x\in\mathbb{R}^n:\operatorname{dist}(x,K)<r\bigr\}$, and let $ V:=K_{2\varepsilon}$. Then $\overline{K_{3\varepsilon}}\subset U$. Define the scaled mollifier
	\[
	\psi_\varepsilon(x)
	:=
	\frac{1}{\varepsilon^n}
	\psi\!\left(\frac{x}{\varepsilon}\right)
	\]
	and set $\phi:=\psi_\varepsilon*\mathbf{1}_V$, where $\mathbf{1}_V$ denotes the characteristic function of $V$. Since $V$ is bounded and $\psi_\varepsilon\in C_c^\infty(\mathbb{R}^n)$, we have $\phi\in C_c^\infty(\mathbb{R}^n)$. Moreover, for every $x\in\mathbb{R}^n$,
	\[
	\phi(x)
	=
	\int_{\mathbb{R}^n}
	\psi_\varepsilon(z)\mathbf{1}_V(x-z)\,\mathrm dz.
	\]
	Since $0\leq\mathbf{1}_V\leq1$, $\psi_\varepsilon\geq0$, and $\int_{\mathbb{R}^n}\psi_\varepsilon(z)\,dz=1$, it follows that $0\leq\phi(x)\leq1$. 
	This implies that  $0\leq\phi\leq1$ in $\Omega$.
	
	We next show that $\phi\equiv1$ on $K$. Let $x\in K$. If $z\in\operatorname{supp}\psi_\varepsilon$, then $|z|\leq\varepsilon$, and hence
	\[
	\operatorname{dist}(x-z,K)
	\leq |z|
	\leq\varepsilon
	<2\varepsilon.
	\]
	Therefore $x-z\in V$, and so $\mathbf{1}_V(x-z)=1$. It follows that 
	\[
	\phi(x)
	=	
	\int_{\mathbb{R}^n}
	\psi_\varepsilon(z)\mathbf{1}_V(x-z)\,\mathrm dz
	=
	\int_{\mathbb{R}^n}\psi_\varepsilon(z)\,\mathrm dz
	=1.
	\]
	Hence $\phi\equiv1$ on $K$. Next, if $\phi(x)\neq0$, then there exists $z\in\mathbb{R}^n$ such
	that $\psi_\varepsilon(z)>0$ and $x-z\in V$. Consequently,
	\[
	\operatorname{dist}(x,K)
	\leq
	\operatorname{dist}(x-z,K)+|z|
	<3\varepsilon.
	\]
	Hence $\operatorname{supp}\phi\subset\overline{K_{3\varepsilon}}\subset U$. Thus, $\phi\in C_c^\infty(\Omega)$ and $\phi=0$ in $\Omega\setminus U$. Finally, by scaling the inequality for $\psi$, we obtain 	
	\[
	\Delta\psi_\varepsilon(x)
	=
	\varepsilon^{-n-2}
	(\Delta\psi)\!\left(\frac{x}{\varepsilon}\right)
	\geq
	-\frac{C_0}{\varepsilon^2}\psi_\varepsilon(x).
	\] 
	Since $\mathbf{1}_V\geq0$, convolution with $\mathbf{1}_V$ preserves this inequality. Therefore, we have
	\[
	\Delta\phi=(\Delta\psi_\varepsilon)*\mathbf{1}_V\geq-\frac{C_0}{\varepsilon^2}
	\psi_\varepsilon*\mathbf{1}_V=-\frac{C_0}{\varepsilon^2}\phi.
	\]
	Taking $a:=\frac{C_0}{\varepsilon^2}>0$, we then  completes the proof.
\end{proof}

\section{Asymptotically vanishing diffusion}
Let \(\Omega\subset\R^n\) be a bounded domain with smooth boundary, $d$ be a positive real number, and \(\eps\in C(\Rp,(0,\infty))\) satisfy
$\lim_{\hat t\to\infty}\eps(\hat t)=0$.   We consider the following impulsive reaction--diffusion equation:
\begin{equation}
\left\{
\begin{aligned}
\frac{\partial u_m}{\partial t}
&=d\eps(mT+t)\Delta u_m+F(t,u_m),
&&0<t\le T,\ x\in\Omega,\\
u_m(0,x)&=g(N_m(x)),
&&x\in\Omega,\\
N_{m+1}(x)&=u_m(T,x),
&&x\in\Omega,\\
Bu_m(t,x)&=0,
&&0<t\le T,\ x\in\partial\Omega.
\end{aligned}
\right.
\label{eq:vanishing-diffusion-system}
\end{equation}
 Throughout this section, we assume that the functions \(F\) and \(g\) satisfy the analogues of (A1)--(A3), and that \(g\) is nondecreasing on \(\mathbb{R}_+\). The boundary operator \(B\) is either the homogeneous Dirichlet operator or the homogeneous Neumann operator.  Let   $Y_+=C(\overline{\Omega}, \mathbb{R}_+)$  in the case of $Bu=\frac{\partial u}{\partial \nu}$,  and  $Y_+=C_0(\overline{\Omega}, \mathbb{R}_+)$  in the case of  $Bu=u$.

Let \(W(t)\) be the solution map of the periodic ordinary differential equation:
\begin{equation}
w'=F(t,w).
\label{eq:limit-ode-F}
\end{equation}
Then \(W:=W(T)\) is the Poincar\'e map of  \eqref{eq:limit-ode-F}, and the limiting impulsive recursion 
associated with \eqref{eq:vanishing-diffusion-system} is
\begin{equation}
M_{m+1}=W\big(g(M_m)\big)=:\widehat W(M_m).
\label{eq:limit-map-W}
\end{equation}
Set
\begin{equation*}
\cR:=g'(0)\exp\left\{\int_0^T F_u(s,0)\,ds\right\}.
\end{equation*}
By Lemma~\ref{lem:kinetic} with \(k=0\) and $f$  replaced by \(F\),  it follows that if \(\cR>1\), then the map \(\widehat W\) admits a unique positive fixed point \(M^*>0\), and the corresponding positive impulsive periodic solution is
\begin{equation*}
w^*(t):=W(t)\big(g(M^*)\big), \qquad 0\le t\le T.
\end{equation*}

For any \(\eta>0\), let
\begin{equation*}
	\Omega_\eta:=\{x\in\Omega:\dist(x,\partial\Omega)>\eta\}.
\end{equation*}
By Lemma \ref {lem:cutoff} with $U=\Omega_\eta$ and $K=\overline{\Omega_{2\eta}}$, it follows that  for every sufficiently small \(\eta>0\), 
there exist \(\phi_\eta\) and \(a_\eta>0\) such that \(\phi_\eta\equiv1\) on \(\Omega_{2\eta}\), \(\phi_\eta\equiv0\) in \(\Omega\setminus\Omega_\eta\), and \(\Delta\phi_\eta\ge-a_\eta\phi_\eta\) in $\Omega$.

For any given small \(\eta>0\), we consider the nonautonomous impulsive kinetic system:
\begin{equation}
\left\{
\begin{aligned}
\frac{dw_m}{dt}
&=F(t,w_m)-a_\eta d\eps(mT+t)w_m,
&&0<t\le T,\\
w_m(0)&=g(M_m),\\
M_{m+1}&=w_m(T).
\end{aligned}
\right.
\label{eq:auxiliary-kinetic}
\end{equation}

\begin{lemma}\label{lem:cutoffsub}
Let \((w_m(t),M_m)\) be a nonnegative solution of \eqref{eq:auxiliary-kinetic}. Then
$(\underline u_m(t,x), \underline N_m(x)):=\phi_\eta(x)(w_m(t), M_m)$
is a subsolution of system \eqref{eq:vanishing-diffusion-system}.
\end{lemma}

\begin{proof}
Since \(F(t,u)/u\) is decreasing in \(u>0\), $F(t,u)$ is subhomogeneous on $\mathbb{R}_+$ in the sense that 
\begin{equation*}
F(t,\alpha u)\ge\alpha F(t,u),
\qquad \forall  \, 0\le\alpha\le1,\quad u\ge0.
\end{equation*}
Thus, we have  
\begin{align*}
&\frac{\partial\underline u_m}{\partial t}
-d\eps(mT+t)\Delta\underline u_m
-F(t,\underline u_m)\\
&\le
\phi_\eta w_m'
-d\eps(mT+t)(\Delta\phi_\eta)w_m
-\phi_\eta F(t,w_m)\\
&\le
\phi_\eta\big[w_m'
+a_\eta d\eps(mT+t)w_m
-F(t,w_m)\big]\\
&=0,\quad  \forall  \, 0<t\le T.
\end{align*}
It remains to verify the impulsive inequality. Since \(g(u)/u\) is nonincreasing in $u>0$, there holds
\[
g(\alpha u)\ge\alpha g(u),
\qquad 0\le\alpha\le1,\quad u\ge0.
\]
It then follows that 
\[
\underline u_m(0,x)=\phi_\eta(x)g(M_m)\le g\big(\phi_\eta(x)M_m\big)=g(\underline N_m(x)).
\]
Moreover,
\[
\underline N_{m+1}(x)
=\phi_\eta(x)M_{m+1}
=\phi_\eta(x)w_m(T)
=\underline u_m(T,x).
\]
Since \(\phi_\eta\) has compact support in \(\Omega\), the boundary condition is also satisfied. 
\end{proof}

Regarding the asymptotic behavior of solutions to the auxiliary system \eqref{eq:auxiliary-kinetic}, we have the following
result.

\begin{lemma}\label{lem:auxconv}
Assume \(\cR>1\). Then every positive solution \((w_m(t),M_m)\) of \eqref{eq:auxiliary-kinetic} satisfies
$\lim_{m\to\infty}M_m=M^*$ and $\lim_{m\to\infty}\sup_{t\in[0,T]}|w_m(t)-w^*(t)|=0$.
\end{lemma}

\begin{proof}
For each $m\geq 0$,  let \(Q_m^\eta\) denote the one-cycle map associated with \eqref{eq:auxiliary-kinetic} such that \(M_{m+1}=Q_m^\eta(M_m)\). Define
\[
\mathcal T_0^\eta=I,\qquad
\mathcal T_m^\eta
=Q_{m-1}^\eta\circ Q_{m-2}^\eta\circ\cdots\circ Q_1^\eta\circ Q_0^\eta,
\qquad \forall \,m\ge1.
\]
Then \(M_m=\mathcal T_m^\eta(M_0)\).  Using  $\lim_{\hat t\to \infty}\eps(\hat t)=0$, together with the continuous dependence of solutions on parameters and initial values, we can verify that 
\[
\lim_{m \to \infty}Q_m^\eta(s)= \widehat W(s)
\]
uniformly for \(s\) in bounded intervals. Thus, \(\{\mathcal T_m^\eta\}_{m\ge0}\) is an asymptotically autonomous discrete process with limit map \(\widehat W\), as defined in \cite{Zhao2017}.  Since  the orbit \(\{M_m:m\ge0\}\) is bounded,  its omega-limit set
\[
\omega(M_0):=\{\xi\in\Rp:\text{there exists }m_j\to\infty
\text{ such that }M_{m_j}\to\xi\}
\]
is nonempty and compact. By \cite[Lemma 1.2.2]{Zhao2017}, \(\omega(M_0)\) is internally chain transitive for \(\widehat W\).

We now claim that \(\omega(M_0)\ne\{0\}\). Since \(\cR>1\), choose \(\theta>0\) sufficiently small such that
\[
q:=\bigl(g'(0)-\theta\bigr)
\exp\left\{\int_0^T\bigl(F_u(s,0)-2\theta\bigr)\,ds\right\}>1.
\]
By the differentiability of $F(t,s)$ and $g(s)$ at $s=0$, there exists \(r_0>0\) such that
\[
F(t,s)\ge\bigl(F_u(t,0)-\theta\bigr)s,
\qquad
g(s)\ge\bigl(g'(0)-\theta\bigr)s, \quad \forall t\in[0,T],\,   s\in [0,r_0].
\]
Since $\lim_{\hat t\to \infty}\eps(\hat t)=0$, there exists an integer \(m_0\) such that
\[
a_\eta d\eps(mT+t)\le\theta,
\qquad \forall \,   m\ge m_0,\, t\in[0,T].
\]
Assume, by contradiction, that \(\omega(M_0)=\{0\}\).  Then  $\lim_{m \to \infty}M_m=0$, and hence, 
the  continuous dependence of soutions on initial values implies  $\lim_{m \to \infty}w_m(t)=0$ uniformly for \(t\in[0,T]\).  
Thus,  we have \(0\le w_m(t)\le r_0\) on \([0,T]\) for all sufficiently large \(m\). The comparison theorem therefore gives
\[
M_{m+1}=w_m(T)
\ge \bigl(g'(0)-\theta\bigr)
\exp\left\{\int_0^T\bigl(F_u(s,0)-2\theta\bigr)\,ds\right\}M_m
=qM_m,
\]
for all sufficiently large \(m\).  Since  \(q>1\),  it follows that  $\lim_{m \to \infty}M_m=\infty$, 
which contradicts  $\lim_{m \to \infty}M_m=0$. This proves the claim.

By Lemma~\ref{lem:kinetic}, applied with \(F\) in place of \(f\) and \(k=0\), \(M^*\) is globally asymptotically stable for \(\widehat W\) in \(\Rp\setminus\{0\}\). Hence \(W^s(M^*)=\Rp\setminus\{0\}\), and the claim implies \(\omega(M_0)\cap W^s(M^*)\ne\emptyset\). Since \(\omega(M_0)\) is internally chain transitive for \(\widehat W\), \cite[Theorem 1.2.1]{Zhao2017} yields \(\omega(M_0)=\{M^*\}\). 
This implies $\lim_{m\to\infty}M_m=M^*$.

It is easy to see  that 
\[
w_m(t)=\Psi_m^\eta\bigl(t,g(M_m)\bigr),\quad \forall  m\geq 0, \,  \, 0\le t\le T,
\]
where \(\Psi_m^\eta(t,\cdot)\) is the time-\(t\) map of the differential equation in \eqref{eq:auxiliary-kinetic}. Since 
$\lim_{m\to\infty}a_\eta d\eps(mT+t)=0$ uniformly for \(t\in[0,T]\) and $\lim_{m\to\infty}M_m=M^*$, it follows 
from the continuous dependence of solutions on parameters and initial values that 
$\lim_{m\to\infty}\Psi_m^\eta\bigl(t,g(M_m)\bigr)=W(t)\bigl(g(M^*)\bigr)$ uniformly for $t\in [0.T]$,
and hence, $\lim_{m\to\infty}\sup_{t\in[0,T]}|w_m(t)-w^*(t)|=0$.
\end{proof}

We are now in a position to prove the main result of this section.

\begin{theorem}\label{thm:vanish}
Assume  that \(F\) and \(g\) satisfy the analogues of \emph{(A1)--(A3)}, and that \(g\) is nondecreasing on \(\mathbb{R}_+\).
Then the following statements are valid:
\begin{enumerate}[label=\rm(\roman*)]
\item If \(\cR\le1\), then every solution of \eqref{eq:vanishing-diffusion-system} with $N_0\in Y_+$
satisfies
$\lim_{m\to\infty}\norm{u_m(t,\cdot)}_\infty=0$ uniformly for \(t\in[0,T]\).
\item If \(\cR>1\), then every solution of \eqref{eq:vanishing-diffusion-system} with $N_0\in Y_+\setminus \{0\}$
satisfies
$\lim_{m\to\infty}u_m(t,x)=w^*(t)$
uniformly for \(t\in[0,T]\) and uniformly for \(x\) in any compact subset of \(\Omega\).
\end{enumerate}
\end{theorem}

\begin{proof}
Let \((u_m,N_m)\) be a given nonnegative solution of \eqref{eq:vanishing-diffusion-system}.
In the case where \(\cR\le1\), we choose \(M_0\ge\norm{N_0}_\infty\) and let \((w_m(t),M_m)\), $m\geq 0$, be the solution of 
the following impulsive kinetic system:
\begin{equation}
	\left\{
	\begin{aligned}
		&\frac{\partial w_m}{\partial t}=F(t,w_m),\quad 0<t\le T,\\
		&w_m(0)=g(M_m),\\
		&M_{m+1}=w_m(T).
	\end{aligned}
	\right.
	\label{ODEsystem}
\end{equation}
Clearly, \(w_m(t)\) is a supersolution of \eqref{eq:vanishing-diffusion-system}. Then the comparison principle gives rise to
\[
0\le u_m(t,x)\le w_m(t), \quad \forall \,  m\ge0, \, t\in[0,T],\,  x\in\Omega. 
\]
Since zero is globally attractive for system \eqref{ODEsystem}, assertion (i) follows.

In the case where \(\cR>1\),  for any given compact subset $K$ of $\Omega$, we choose a small number \(\eta>0\) such that \(K\subset\Omega_{2\eta}\).   Let $N_0\in Y_+\setminus \{0\}$ be given. By  
the strong maximum principle (see, e.g., \cite[Theorem 7.4.1]{Smith1995}), it then follows that 
for every  $m\geq 1$, $N_{m}(\cdot)\in {\rm Int}(C(\overline{\Omega},\mathbb{R}_+))$ 
in the case of $Bu=\frac{\partial u}{\partial \nu}$, and $N_{m}(\cdot)\in {\rm Int}(C_0^1(\overline{\Omega},\mathbb{R}_+))$ 
in the case of  $Bu=u$. Thus,  we can choose a real number \(M_{1}>0\) sufficiently small such that
\begin{equation*}
\phi_\eta(x)M_{1}\le N_{1}(x),
\quad \forall \, x\in\Omega.
\end{equation*}
Let \((w_m,M_m)\) be the solution of \eqref{eq:auxiliary-kinetic} starting from $m=1$. By Lemma~\ref{lem:cutoffsub} and the comparison principle,
it then follows that 
\[
u_m(t,x)\ge\phi_\eta(x)w_m(t),
\quad \forall \, m\ge 1,  \, t\in [0,T],\,  x\in\Omega.
\]
Since \(\phi_\eta\equiv1\) on \(K\), Lemma~\ref{lem:auxconv} yields
\begin{equation}
\liminf_{m\to\infty}\inf_{x\in K}
\big[u_m(t,x)-w^*(t)\big]\ge0
\label{eq:lower-limit-vanishing}
\end{equation}
uniformly for \(t\in[0,T]\).

For the upper estimate, choose \(\overline M_{1}\ge\norm{N_{1}}_\infty\), and let \((\overline w_m,\overline M_m)\) 
be the solution of system \eqref{ODEsystem} for all $m\geq 1$. Then \(\overline w_m(t)\) is a supersolution of \eqref{eq:vanishing-diffusion-system}, 
and hence,
\[
u_m(t,x)\le\overline w_m(t), \quad \forall m\ge 1, \, t\in [0,T],\,  x\in\Omega.
\]
By the global attractivity of \(w^*(t)\) for system \eqref{ODEsystem}, we obtain
\begin{equation}
\limsup_{m\to\infty}\sup_{x\in\Omega}
\big[u_m(t,x)-w^*(t)\big]\le0
\label{eq:upper-limit-vanishing}
\end{equation}
uniformly for \(t\in[0,T]\). Combining \eqref{eq:lower-limit-vanishing} and \eqref{eq:upper-limit-vanishing} proves 
assertion (ii).
\end{proof}

\begin{remark}
Under the homogeneous Dirichlet boundary condition, the convergence in Theorem~\ref{thm:vanish}(ii) is locally uniform in $\Omega$, because the positive solution of the limiting kinetic system  \eqref{ODEsystem} 
does not satisfy the zero boundary condition.  For the homogeneous Neumann condition,  
we can choose  any solution  of  \eqref{ODEsystem}  as a spatially homogeneous 
subsolution of system \eqref{eq:vanishing-diffusion-system}  to obtain the uniform
convergence on $\overline{\Omega}$.
\end{remark}

\section{The full impulsive system}
We first study the global dynamics of system \eqref{ReduceE} in the case where \(g\) is nondecreasing on $\mathbb{R}_+$.
Consider the spatially homogeneous nonautonomous impulsive system:
\begin{equation}
\left\{
\begin{aligned}
\frac{dz_m}{dt}
&=f(t,z_m)-n\frac{\rho'(mT+t)}{\rho(mT+t)}z_m,
&&0<t\le T,\\
z_m(0)&=g(Z_m),\\
Z_{m+1}&=z_m(T).
\end{aligned}
\right.
\label{eq:full-kinetic}
\end{equation}
Let \(P_m\) be its one-cycle map such that \(Z_{m+1}=P_m(Z_m)\). Define
\begin{equation*}
\mathcal S_0=I,\qquad
\mathcal S_m=P_{m-1}\circ P_{m-2}\circ\cdots\circ P_1\circ P_0,
\qquad \forall  \,  m\ge1.
\end{equation*}
It then follows that 
\begin{equation*}
Z_m=\mathcal S_m(Z_0),\qquad \forall  \, m\ge0.
\end{equation*}
In view of (A5), we can verify that $\lim_{m\to\infty}P_m(s)=\widehat W_k(s)$ uniformly for \(s\) in bounded intervals. Therefore, \(\{\mathcal S_m\}_{m\ge0}\) is an asymptotically
autonomous discrete process with limit map \(\widehat W_k\).

\begin{lemma}\label{lem:odefull}
Assume that \emph{(A1)--(A3)} and \emph{(A5)} hold, and \(g\) is nondecreasing on $\mathbb{R}_+$. Then the following statements are valid:
\begin{enumerate}[label=\rm(\roman*)]
	\item If \(\cR_k\le1\), then every nonnegative solution of \eqref{eq:full-kinetic} satisfies
	$\lim_{m\to\infty}z_m(t)=0$ uniformly for \(t\in[0,T]\).
	\item If \(\cR_k>1\), then every positive solution of \eqref{eq:full-kinetic} satisfies
	$\lim_{m\to\infty}z_m(t)=w_k^*(t)$ uniformly for \(t\in[0,T]\).
\end{enumerate}
\end{lemma}

\begin{proof}
In view of (A1)--(A3) and Lemma \ref{lem:bounded}, we see that for any
\(Z_0\ge0\), the forward orbit of \(\{\mathcal S_m\}_{m\ge0}\)
is bounded and its omega-limit set $\omega(Z_0)$ is nonempty, compact
and internally chain transitive for the map \(\widehat W_k\) by \cite[Lemma 1.2.2]{Zhao2017}.
In the case where \(\cR_k\le1\), Lemma~\ref{lem:kinetic} (i) shows that zero is globally asymptotically stable for
\(\widehat W_k\), and hence,  \(W^s(0)=\Rp\). Clearly, \(\omega(Z_0)\cap W^s(0)\ne\emptyset\). By \cite[Theorem 1.2.1]{Zhao2017}, we have 
\(\omega(Z_0)=\{0\}\).  Thus, $\lim_{m\to \infty}Z_m=0$.  It then follows that 
$\lim_{m\to\infty}z_m(t)=0$ uniformly for \(t\in[0,T]\).
In the case where \(\cR_k>1\), we can use arguments similar to those in Lemma \ref{lem:auxconv}
to prove that  for any $Z_0>0$, we have $\lim_{m\to\infty}Z_m=M_k^*$ 
and $\lim_{m\to\infty}\sup_{0\le t\le T}|z_m(t)-w_k^*(t)|=0$.
\end{proof}

\begin{theorem}\label{thm:mainmono}
Assume that   \emph{(A1)--(A3)} and \emph{(A5)} hold, and \(g\) is nondecreasing on $\mathbb{R}_+$. 
Let \(u_m(t,y,\phi)\) be the solution of system \eqref{ReduceE} with initial value $N_0=\phi\in X_+$. Then the following statements are valid:
\begin{enumerate}[label=\rm(\roman*)]
\item If \(\cR_k\le1\), then  for every $\phi\in X_+$,
$\lim_{m\to\infty}\norm{u_m(t,\cdot,\phi)}_\infty=0$ uniformly for \(t\in[0,T]\).
\item If \(\cR_k>1\), then for every $\phi\in X_+\setminus \{0\}$,
$\lim_{m\to\infty}u_m(t,y,\phi)=w_k^*(t)$ uniformly for \(t\in[0,T]\) and uniformly for \(y\) in any compact 
subset of \(\Omega_0\).
\end{enumerate}
\end{theorem}

\begin{proof}
Let \((u_m(t,y),N_m(y))\) be the unique solution of \eqref{ReduceE} with $N_0=\phi$.
In the case where \(\cR_k\le1\), we choose a solution of system \eqref{eq:full-kinetic} whose initial value 
$Z_0$ dominates \(N_0\). Since any spatially homogeneous solution of \eqref{eq:full-kinetic} is a supersolution of 
system \eqref{ReduceE}, the comparison principle and Lemma~\ref{lem:odefull}(i) imply $\lim_{m\to\infty}\norm{u_m(t,\cdot)}_\infty=0$ uniformly for \(t\in[0,T]\).

In the case where \(\cR_k>1\), we fix a real number \(\delta>0\) sufficiently small such that
\begin{equation*}
\cR_{k+\delta}
=g'(0)\exp\left\{\int_0^T\big(f_u(s,0)-n(k+\delta)\big)\,ds\right\}>1.
\end{equation*}
In view of  (A5), there exists an integer \(m_\delta\) such that
\begin{equation*}
\frac{\rho'(mT+t)}{\rho(mT+t)}<k+\delta,
\qquad \forall \, m\ge m_\delta,\quad t\in[0,T].
\end{equation*}
It then follows that  for all   \(m\ge m_\delta\),  \(u_m(t,y)\) is a supersolution of the following impulsive 
system:
\begin{equation}
\left\{
\begin{aligned}
\frac{\partial U_m}{\partial t}
&=\frac{d}{\rho^2(mT+t)}\Delta U_m
+f(t,U_m)-n(k+\delta)U_m,
\qquad 0<t\le T,\\
U_m(0,y)&=g(L_m(y)),\\
L_{m+1}(y)&=U_m(T,y),\\
BU_m&=0.
\end{aligned}
\right.
\label{eq:comparison-system-delta}
\end{equation}
Since $\lim_{\hat t \to \infty}\frac{1}{\rho^2(\hat t)}=0$, Theorem~\ref{thm:vanish}, as applied to sysyem \eqref{eq:comparison-system-delta},  implies that for every compact subset $K$ of $\Omega_0$,
\begin{equation}
\liminf_{m\to\infty}\inf_{y\in K}
\big[u_m(t,y)-w_{k+\delta}^*(t)\big]\ge0
\label{eq:lower-limit-full}
\end{equation}
uniformly for \(t\in[0,T]\). For the upper estimate, the  solution of system \eqref{eq:full-kinetic} 
with a sufficiently large initial value $Z_0$ is a supersolution of system \eqref{ReduceE}. Lemma~\ref{lem:odefull}(ii) then yields
\begin{equation}
\limsup_{m\to\infty}\sup_{y\in\Omega_0}
\big[u_m(t,y)-w_k^*(t)\big]\le0
\label{eq:upper-limit-full}
\end{equation}
uniformly for \(t\in[0,T]\).
Since the fixed point of \(\widehat W_{k+\delta}\) and the associated periodic solution 
depend
 continuously on \(\delta\) for small \(\delta\), it follows that $\lim_{\delta\to0^+}w_{k+\delta}^*(t)=w_k^*(t)$ uniformly for \(t\in[0,T]\). Letting \(\delta\to0^+\) in \eqref{eq:lower-limit-full} and combining with \eqref{eq:upper-limit-full}, we obtain 
$\lim_{m\to\infty}u_m(t,y)=w_k^*(t)$ uniformly for \(t\in[0,T]\) and uniformly for \(y\) in $K$.
\end{proof}

\begin{remark}
In view of $\cR_k=\cR_0e^{-nkT}$, it follows  that 
when the other parameters are fixed, an increase in the asymptotic relative expansion rate \(k\) decreases the threshold 
value because of the dilution effect. In the special case \(k=0\), the threshold coincides with that of the spatially homogeneous impulsive model without dilution.
\end{remark}

Next, we consider the case where \(g\) is  not monotone on $\mathbb{R}_+$. Define
\begin{equation*}
g^+(N):=\max_{0\le v\le N}g(v),
\qquad N\ge0.
\end{equation*}
Then \(g^+\) is nondecreasing and locally Lipschitz, \(g^+(N)/N\) is nonincreasing, and (A4) implies 
$(g^+)'(0)=g'(0)$.
Replacing \(g\) with \(g^+\) in the limiting kinetic system \eqref{eq:limiting-kinetic}, we see from  Lemma~\ref{lem:kinetic}  that if \(\cR_k>1\), the resulting scalar recursion admits a unique positive fixed point \(\gamma_k^+\). We further define
\begin{equation*}
g^-(N):=\min_{N\le v\le\gamma_k^+}g(v),
\qquad 0\le N\le\gamma_k^+.
\end{equation*}
Again, \(g^-\) is nondecreasing and locally Lipschitz, \(g^-(N)/N\) is nonincreasing, and
$(g^-)'(0)=g'(0)$.
Moreover, there holds
\begin{equation*}
g^-(N)\le g(N)\le g^+(N), \qquad 0\le N\le\gamma_k^+.
\end{equation*}

Let \(w_k^{+*}(t)\) and \(w_k^{-*}(t)\) be the positive periodic impulsive solutions of the limiting kinetic systems associated with \(g^+\) and \(g^-\), respectively.  Then we have the following result.

\begin{theorem}\label{thm:nonmono}
Assume \emph{(A1)--(A5)} hold. Then the following statements are valid for system \eqref{ReduceE}:
\begin{enumerate}[label=\rm(\roman*)]
\item If \(\cR_k\le1\), then every nonnegative solution satisfies
$\lim_{m\to\infty}\norm{u_m(t,\cdot)}_\infty=0$ uniformly for \(t\in[0,T]\).
\item If \(\cR_k>1\), then for every nontrivial nonnegative solution and every compact subset $K$ of  $\Omega_0$,
$\limsup_{m\to\infty}\sup_{(t,y)\in [0,T]\times K}\big[u_m(t,y)-w_k^{+*}(t)\big]\le0$ 
and
$\liminf_{m\to\infty}\inf_{(t,y)\in [0,T]\times K}\big[u_m(t,y)-w_k^{-*}(t)\big]\geq 0$.
\end{enumerate}
\end{theorem}

\begin{proof}
In the case where \(\cR_k\le1\), let \(u_m^+(t,y)\) be the solution of the upper comparison system obtained by replacing \(g\) with \(g^+\) and taking the same initial function. Then the comparison principle implies 
\[
0\le u_m(t,y)\le u_m^+(t,y),\qquad \forall m\geq 0,\,  t\in [0,T], \, y\in \Omega_0.
\]
Since the upper system has the same threshold value \(\cR_k\), Theorem~\ref{thm:mainmono}(i) gives extinction of \(u_m^+\), and hence of \(u_m\).

In the case where \(\cR_k>1\), let \(u_m^+(t,y)\) and \(u_m^-(t,y)\) be the solutions of 
the impulsive systems generated by \(g^+\) and \(g^-\), respectively. Choose the lower initial function sufficiently small so that it lies below the original initial profile and takes values in \([0,\gamma_k^+]\). Starting the lower comparison system at a sufficiently large cycle if necessary, the comparison principle implies that 
\[
0\le u_m^-(t,y)\le u_m(t,y)\le u_m^+(t,y),\qquad \forall \, t\in[0,T],\, y\in \Omega_0,
\]
for all subsequent cycles. Applying Theorem~\ref{thm:mainmono}(ii) to these two monotone systems, we obtain the desired result.
\end{proof}

\section{Discussion}

In this paper, we have studied a class of  impulsive reaction--diffusion population models with  an asymptotically unbounded domain. Earlier studies have treated impulsive models on particular growing domains, including the logistic harvesting model in \cite{ZhangZhu2025}, as well as general impulsive systems on asymptotically bounded and asymptotically periodic domains \cite{ZhuZhao2025,FanZhao2026}. The contribution here concerns the general asymptotically unbounded regime. After transforming the model to a fixed reference domain, the diffusion coefficient tends to zero and the dilution coefficient converges to \(nk\). Accordingly, the evolution  dynamics is governed by a spatially homogeneous impulsive kinetic system rather than by a limiting parabolic equation.

The major analytical issue arises in the persistence case. For nonimpulsive systems with vanishing diffusion, the lower estimate in \cite{LamZhaoZhu2024} was obtained from a generalized subsolution involving a connecting orbit of the limiting kinetic system. Here we construct  a new  subsolution by combining a cut-off function \(\phi_\eta\) satisfying \(\Delta\phi_\eta\ge-a_\eta\phi_\eta\) with an auxiliary impulsive kinetic system containing the vanishing correction \(-a_\eta d\eps(t)w\). The subhomogeneity of both the reaction term and the birth function then guarantees the differential and impulsive inequalities. This  enable us to prove  the  local uniform convergence of solutions of the impulsive reaction-diffusion system 
to the positive periodic solution of its  limiting impulsive kinetic system.

Our threshold value
\[
\cR_k
=g'(0)\exp\left\{\int_0^T f_u(s,0)\,ds-nkT\right\}
\]
separates the effects of the birth pulse, within-cycle growth, and habitat expansion. The asymptotic diffusion coefficient vanishes and therefore does not enter this threshold explicitly. Nevertheless, diffusion continues to influence the transient spatial profile and, under the Dirichlet boundary condition, produces a boundary layer which explains why persistence convergence is naturally local in the interior of the reference domain.

It would be interesting to extend the present approach to spatially heterogeneous reaction and birth terms, to cooperative multi-species impulsive systems, and to more general anisotropic domain deformations. Our
constructed  subsolution here may also be useful for other hybrid parabolic equations whose diffusion coefficients vanish asymptotically.

\

\noindent
{\bf Acknowledgments.} Zhao's research was supported in part by the NSERC of Canada (RGPIN-2025-04963).


\begin{thebibliography}{99}
\bibitem{BaiLouZhao2022}
Z.~Bai, Y.~Lou, and X.-Q.~Zhao,
Spatial dynamics of species with annually synchronized emergence of adults,
\emph{J. Nonlinear Sci.} 32 (2022), 78.

\bibitem{FanZhao2026}
S.~Fan and X.-Q.~Zhao,
Threshold dynamics of impulsive reaction--diffusion systems with asymptotically periodic domain,
\emph{Z. Angew. Math. Phys.} 77 (2026), 259.

\bibitem{FazlyLewisWang2017}
M.~Fazly, M.~A.~Lewis, and H.~Wang,
On impulsive reaction--diffusion models in higher dimensions,
\emph{SIAM J. Appl. Math.} 77 (2017), 224--246.

\bibitem{FazlyLewisWang2020}
M.~Fazly, M.~A.~Lewis, and H.~Wang,
Analysis of propagation for impulsive reaction--diffusion models,
\emph{SIAM J. Appl. Math.} 80 (2020), 521--542.

%\bibitem{Zhao2001}
%M.W. Hirsch, H.L. Smith and X.-Q. Zhao, Chain transitivity, attractivity
%and strong repellors for semidynamical systems,  \emph{J. Dyn. Diﬀer. Equ.}
 %13 (2001), 107–-131.

\bibitem{LamZhaoZhu2024}
K.-Y.~Lam, X.-Q.~Zhao and M.~Zhu,
Global dynamics of reaction--diffusion systems with a time-varying domain,
\emph{SIAM J. Appl. Math.} 84 (2024), 1742--1765.

\bibitem{LewisLi2012}
M.~A.~Lewis and B.~Li,
Spreading speed, traveling waves, and minimal domain size in impulsive reaction--diffusion models,
\emph{Bull. Math. Biol.} 74 (2012), 2383--2402.

\bibitem{MengLinPedersen2021}
Y.~Meng, Z.~Lin, and M.~Pedersen,
Effects of impulsive harvesting and an evolving domain in a diffusive logistic model,
\emph{Nonlinearity} 34 (2021), 7005--7029.

\bibitem{Smith1995}
H.L. Smith, \emph{Monotone Dynamical Systems, An Introduction to the 
Theory of Competitive and Cooperative Systems}, American Mathematical Society, 
1995.

\bibitem{WuZhao2019}
R.~Wu and X.-Q.~Zhao,
Spatial invasion of a birth pulse population with nonlocal dispersal,
\emph{SIAM J. Appl. Math.} 79 (2019), 1075--1097.

\bibitem{WuZhao2022}
R.~Wu and X.-Q.~Zhao,
The evolution dynamics of an impulsive hybrid population model with spatial heterogeneity,
\emph{Commun. Nonlinear Sci. Numer. Simul.} 107 (2022), 106181.

\bibitem{ZhangYiChen2024}
Y.~Zhang, T.~Yi, and Y.~Chen,
Spreading dynamics of an impulsive reaction--diffusion model with shifting environments,
\emph{J. Differential Equations} 381 (2024), 1--19.

\bibitem{ZhangZhu2025}
H.~Zhang and M.~Zhu,
Population dynamics of a logistic model incorporating harvesting pulses on a growing domain,
\emph{Commun. Nonlinear Sci. Numer. Simul.} 146 (2025), 108768.

\bibitem{Zhao2017}
X.-Q.~Zhao,
\emph{Dynamical Systems in Population Biology}, 2nd ed., Springer,  2017.

%\bibitem{Lin2025}
%Q. Zhou, Z. Lin  and C. A. Santos, On an impulsive faecal-oral model in a periodically evolving environment,
%\emph{Chaos, Solitons $\&$ Fractals} 191(2025), 115825.

\bibitem{ZhuZhao2025}
M.~Zhu and X.-Q.~Zhao,
An impulsive reaction--diffusion model with asymptotically bounded domain,
\emph{J. Math. Anal. Appl.} 552 (2025), 129723.

\end{thebibliography}
\end{document}